\documentclass[12pt,a4paper]{amsart} 
\usepackage{latexsym,amssymb,exscale}
\renewcommand{\ss}{\scriptscriptstyle} 
\renewcommand{\d}{\displaystyle} 
\renewcommand{\ni}{\noindent} 
\renewcommand{\epsilon}{\varepsilon}

\renewcommand{\leq}{\leqslant}

\renewcommand{\geq}{\geqslant}
\newcommand{\V}{{\boldsymbol{V}}} 
 
\renewcommand{\L}{{\mathcal{L}}}

\newcommand{\1}{{\boldsymbol{1}}} 
\newcommand{\F}{{\mathcal{F}}} 
 
\renewcommand{\S}{{T}} 
\renewcommand{\phi}{{\varphi}}

\newcommand{\cZ}{{\mathcal{Z}}}
\newcommand{\N}{{\mathbb{N}}}

\newcommand{\scirc}{{\scriptstyle\circ}}

 \font\bigmath=cmmi10 scaled \magstep2 
	\newcommand{\bigchi}{\hbox{\bigmath \char31}} 
	\newcommand{\ldelta}[1]{{\L}_\Delta^{#1}} 
	\newcommand{\var}{{\mathrm{var}}}

\newcommand{\zmod}[1]{\,(\mathrm{mod}\;{#1})}

\newtheorem{theorem}{Theorem} 
\newtheorem{corollary}[theorem]{Corollary} 
\newtheorem{lemma}[theorem]{Lemma} 
\newtheorem{proposition}[theorem]{Proposition}

\theoremstyle{definition}

\newtheorem*{example}{Example}

\title[Periodic subsystems of finite type]{On the asymptotic measure of periodic subsystems of finite type in symbolic dynamics}
\author{Jean-Ren\'e Chazottes, Zaqueu Coelho and Pierre Collet} %
\address{Centre de Physique Th\'eorique\\
CNRS UMR 7644\\
Ecole Polytechnique\\
F-91128 Palaiseau Cedex\\
FRANCE} %
\email{jean-rene.chazottes@cpht.polytechnique.fr} %
\address{Department of Mathematics \\ University of York\\ Heslington \\
York YO10 5DD\\ UK} %
\email{zc3@york.ac.uk} %
\address{Centre de Physique Th\'eorique\\
CNRS UMR 7644\\
Ecole Polytechnique\\
F-91128 Palaiseau Cedex\\
FRANCE} %
\email{pierre.collet@cpht.polytechnique.fr} %
\date{\today
} %

\begin{document}

\begin{abstract}
Let $\Delta\subsetneq\V$ be a proper subset of the vertices $\V$ of
the defining graph of an irreducible and aperiodic shift of finite type
$(\Sigma_{A}^{+},\S)$.  Let $\Delta_{n}$ be the union of cylinders in
$\Sigma_{A}^{+}$ corresponding to the points $x$ for which the first
$n$-symbols of $x$ belong to $\Delta$ and let $\mu$ be an equilibrium
state of a H\"older potential $\phi$ on $\Sigma_{A}^{+}$.  We know
that $\mu(\Delta_{n})$ converges to zero as $n$ diverges.  We study
the asymptotic behaviour of $\mu(\Delta_{n})$ and compare it with the
pressure of the restriction of $\phi$ to $\Sigma_{\Delta}$.  The
present paper extends some results in~\cite{CCC} to the case when
$\Sigma_{\Delta}$ is irreducible and periodic.  We show an explicit
example where the asymptotic behaviour differs from the aperiodic
case.
\end{abstract}

\maketitle

\setlength{\baselineskip}{15pt}
\enlargethispage*{10pt}
\section*{Introduction}
% intro
\ni The present work was triggered by our study of the convergence
of certain hitting time processes to (marked) Poisson processes \cite{CCC}.
The general setup in the study of asymptotic time distributions
for an ergodic dynamical system $(\Omega,\mu,\S)$ is the following (see e.g. \cite{Coe}). 
Usually one considers a sequence of measurable subsets $B_{n}$ with
asymptotically vanishing measure and one tries to study the asymptotic
behaviour of the random variables:
\[
	\tau_{n}(\omega)\;=\;\inf\{k\geq 1\colon\:\S^{k}(\omega)\in B_{n}\}
\]
when suitably rescaled.
The rescaling is necessary since in general $\tau_{n}(\omega)\to\infty$ for $\mu$-almost every $\omega$.
The most studied case is when $\Omega$ is equipped with a partition and
$B_n=B_n(\omega)$ are cylinder sets about a generic point $\omega$ of $\mu$ (so that
$\cap_n B_n=\omega$). It turns out that $\mu(B_n)$
is the right scaling factor to obtain, for ``sufficiently" mixing systems, a convergence to
a Poisson law for the point process associated with the random variables $\mu(B_n)\tau_n$.
The reader can find a list of references in e.g. \cite{CCC}.
We were interested in the case where the intersection of the $B_n$'s is a non-trivial invariant set.
In \cite{CCC} we consider the case where that intersection is a subsystem of finite type of a shift of finite type.
This is closely related with the so-called Pianigiani-Yorke measures \cite{CMS}.

More precisely, the setup is as follows. 
Let $\Delta\subsetneq\V$ be a proper subset of the vertices $\V$ of
the defining graph of an irreducible and aperiodic shift of finite type
$(\Sigma_{A}^{+},\S)$.  This induces a subsystem of finite type which
we denote by $\Sigma_\Delta$.
Let $\Delta_{n}$ be the union of cylinders in $\Sigma_{A}^{+}$
corresponding to the points $x$ for which the first
$n$-symbols of $x$ belong to $\Delta$.  
Note that $\Sigma_{\Delta}=\cap_{n\geq 0}\,\Delta_{n}$.
Now let $\mu$ be an equilibrium state of a H\"older potential $\phi$
on $\Sigma_{A}^{+}$.
We know that $\mu(\Delta_{n})$ converges to zero as $n$ diverges.
In~\cite{CCC} we assume that $\Sigma_\Delta$ is an irreducible
and {\em aperiodic} subshift of finite type. Then we prove that 
the point process corresponding to the times of hitting $\Delta_{n}$
scaled by $\mu(\Delta_n)$ converges to a (marked) Poisson point process. 
In that paper, we study the asymptotic behaviour of $\mu(\Delta_{n})$ and
compare it with the pressure of the restriction of $\phi$ to $\Sigma_{\Delta}$, call it $P_\Delta$.
We show that $e^{n(P(\phi)-P_\Delta)}\mu(\Delta_n)$ has a limit (which can be identified) as $n\to\infty$.
In the present paper we extend some results in~\cite{CCC} to the case
when $\Sigma_{\Delta}$ is irreducible {\em but periodic}.  This extension
turns out to be non-trivial since, in general, $e^{n(P(\phi)-P_\Delta)}\mu(\Delta_n)$
may {\em not converge} as $n\to\infty$, contrarily to the aperiodic case. Indeed, we provide an explicit example where
this phenomenon appears.

\section{Preliminaries}
Let $\V=\{1,\ldots,\ell\}$ be a finite set of symbols (i.e.~the base
alphabet).  We will assume that $A$ is an aperiodic 
$0$-$1$ $\ell\times\ell$ matrix which defines the allowable 
transitions in a directed graph $\mathcal{G}$ of labelled vertices 
$\V$.  Define the space of one-sided allowable paths in the graph 
$\mathcal{G}$ by
\[
	\Sigma_{A}^{+}\;=\; \big\{x=(x_{n})\in\V^{\N}\colon\; 
	A(x_{i-1},x_{i})=1\,,\;\forall\,i\geq 1\big\}\;.
\]
The space $\Sigma_{A}^{+}$ is compact and metrisable when endowed with
the Tychonov product topology (generated by the discrete topology on
$\V$).  The \emph{shift} $\S$ (of finite type) is the map
$\S\colon\Sigma_{A}^{+}\to\Sigma_{A}^{+}$ defined by
$\S(x)_{n}=x_{n+1}$ for all $n\geq 0$.  This map is continuous and
surjective.  The cylinders, denoted by
\[
	C[i_{0},\ldots,i_{m}]_{k}\;=\; \big\{ x\in\Sigma_{A}^{+}\colon\; 
	x_{j+k}=i_{j}\,,\;\forall\,j=0,\ldots,m \big\}\;,
\]
form a base of open (and closed) sets in $\Sigma_{A}^{+}$.  Let 
$C(\Sigma_{A}^{+})$ denote the space of complex valued continuous 
functions on $\Sigma_{A}^{+}$.  For $\psi\in C(\Sigma_{A}^{+})$, 
consider $\var_{n}(\psi) = \sup\{|\psi(x) - \psi(y)|\colon\: x_{i} = 
y_{i},\: i\leq n\}$.  Given $0 < \theta < 1$, define $|\psi|_{\theta} 
= \sup\big\{\var_{n}(\psi)/\theta^{n}\}$.  The space $\F_{\theta}^{+}= 
\{\psi\in C(\Sigma_{A}^{+})\colon\: |\psi|_{\theta} < \infty\}$ is a 
Banach space when endowed with the norm $\|\psi\|_{\theta} = 
\|\psi\|_{\infty} + |\psi|_{\theta}$, where $\|\cdot\|_{\infty}$ 
denotes the supremum norm.  The union $\F=\cup_{\theta}\,\F_{\theta}$ 
is referred to as the space of H\"older continuous functions on 
$\Sigma_{A}^{+}$.

Let $\phi\in\F_{\theta}^{+}$ be a real valued function.  Define the 
transfer operator $\L_{\phi}$ acting on 
$\F_{\theta}^{+}$ by 
\[
	(\L_{\phi}\psi)(x) = \sum_{\S y=x}\, 
	e^{\phi(y)}\psi(y)\;.
\]
The operator $\L_{\phi}$ has a maximum positive eigenvalue
$e^{P(\phi)}$, which is simple and isolated.  Furthermore, the rest of
the spectrum is contained in a disc of radius strictly less than
$e^{P(\phi)}$ (cf.~\cite{Rue}, \cite{PP}).  The number $P= P(\phi)$ is
called the \emph{pressure} of $\phi$.  There is a unique
$\S$-invariant probability measure $\mu=\mu_{\phi}$ such that
\[
	P(\phi) 
	= h(\mu) + \int \phi\,d\mu\;,
\]	
where $h(\mu)$ denotes the measure-theoretic entropy of $(\S,\mu)$.  
The pressure $P(\phi)$ can also be characterised as the maximum of 
$h(m) + \int\phi\,dm$ over all $\S$-invariant probabilities~$m$.  The 
measure $\mu$ is called the {\em equilibrium state} of $\phi$.  An 
eigenfunction $w$ of $\L_{\phi}$ corresponding to $e^{P(\phi)}$ may be 
taken to be strictly positive, in fact one may take $w$ to be the 
function
\[
	w \;=\; \lim_{n\to\infty}\, e^{-nP(\phi)}\,\L_{\phi}^{n}(\1) \;,
\] 
where $\1$ denotes the constant function equal to $1$.  Replacing 
$\phi$ by $\phi' = \phi - P(\phi) + \log(w) - \log(w\scirc\S)$, we see 
that $\L_{\phi'}\1=\1$ and $P(\phi')=0$.  In this case we say that 
$\phi'$ is \emph{normalised}.  It is easy to see that $\phi$ and 
$\phi'$ have the same equilibrium state $\mu$.  In what follows we 
will assume that $\phi$ is normalised.  Note that in this case the 
transfer operator $\L_{\phi}$ satisfies
\begin{align}
        \nonumber
	\int \psi \: d\mu 
	\;&=\; \int  \L_\phi(\psi) \: d\mu\,, \\
	\int \psi_1\cdot(\psi_2 \scirc \S) \: d\mu \;&=\; \int 
	\L_\phi(\psi_1)\cdot\psi_2  \: d\mu \,,
	\nonumber
\end{align}
for all $\psi,\psi_{1},\psi_{2}\in C(\Sigma_{A}^{+})$.  

Let $\Delta\subseteq\V$ be a sub-alphabet such that $\Delta\not=\V$.  
Consider the closed $\S$-invariant subset $\Sigma_{\Delta}\subseteq 
\Sigma_{A}^{+}$ given by
\begin{equation}
	\Sigma_{\Delta}\;=\; \big\{ x\in\Sigma_{A}^{+}\colon\; 
	x_{i}\in\Delta,\;\forall\,i\geq 0 \big\}\;.
	\nonumber
\end{equation}

In this paper we will consider only the case when $\Sigma_{\Delta}$ is
an irreducible subshift of finite type in its alphabet $\Delta$.  This
means that the restriction of the matrix $A$ to the symbols
of~$\Delta$ defines a matrix $A_{\Delta}$ which is irreducible.  In
particular, the restriction of the shift transformation~$T$ to
$\Sigma_{\Delta}$ is topologically transitive in the induced topology
from $\Sigma_{A}^{+}$.

Let $\phi_{\Delta}$ denote the restriction of $\phi$ to the subsystem 
$\Sigma_{\Delta}$.  Let $P_{\Delta}$ be the pressure of 
$\phi_{\Delta}$ with respect to the subsystem $(\Sigma_{\Delta},\S)$.  
(Note that since $\phi$ is assumed to be normalised we have 
$P(\phi)=0$, therefore $P_{\Delta}<0$.)  Let $\mu_\Delta$ denote the 
equilibrium state of $\phi_{\Delta}$ with respect to the subsystem 
$(\Sigma_{\Delta},{\S})$.  Let $w_{\Delta}$ be the strictly positive 
H\"older continuous function defined on~$\Sigma_{\Delta}$ by
\begin{equation}
	w_{\Delta} \;=\; \lim_{n\to\infty}\, 
	e^{-nP_{\Delta}}\,\L_{\phi_{\Delta}}^{n}(\1) \;.
	\label{eq:wdel}
\end{equation}
Now define the \emph{restricted transfer 
operator} $\ldelta{}$ acting on the space of H\"older continuous 
functions $\F_{\theta}^{+}$ by
\[
	\ldelta{} \psi\;=\;{\L}_\phi(\psi\cdot\bigchi_\Delta)\;,
\]
and consider the subset of~$\Sigma_{A}^{+}$ given by
\begin{equation}
	\cZ_{\Delta} \;=\; \{x\in\Sigma_{A}^{+}\colon\: 
	\exists\,b\in\Delta\,, \;A(b,x_{0})=1 \} \;.
	\label{eq:czdef}
\end{equation}
Note that since $A$ is irreducible and aperiodic in the full 
alphabet~$\V$, $\cZ_{\Delta}$ is a non-empty finite union of cylinder 
sets of~$\Sigma_{\Delta}^{+}$.  In particular, since $\mu$ is fully 
supported on~$\Sigma_{A}^{+}$ we have $\mu(\cZ_{\Delta})>0$.

An improvement to main result of~\cite{CMS}, which is proved
in~\cite{CCC}, gives the following
result.   %
\begin{proposition}
There exists a \emph{unique} H\"older continuous function $h_\Delta$ 
defined on the \emph{whole space} $\Sigma_{A}^{+}$ such that
\[
	\L_{\Delta}(h_{\Delta}) \;=\;  e^{P_{\Delta} } 
	\,h_{\Delta} \;,
\]
and $h_{\Delta}|_{\Sigma_{\Delta}}\equiv w_{\Delta}$, where 
$w_{\Delta}$ is given by~(\ref{eq:wdel}).  The function $h_{\Delta}$ 
is strictly positive on~$\cZ_{\Delta}$ and it is zero on the 
complement~$\cZ_{\Delta}^{c}$.  Moreover,
\[
	\Big\Vert e^{-nP_{\Delta}}\L_{\Delta}^n(\psi) - 
	h_{\Delta} \int_{\Sigma_{\Delta}} \psi \: d\mu_{\Delta} 
	\Big\Vert_{{\Sigma_{A}^{+}}} \; \underset{n \to 
	\infty}{\longrightarrow}\; 0 \;,
\]
for all $\psi \in C({\Sigma_{A}^{+}})$.
\label{prop:cms}
\end{proposition}

The Borel measure $\mu_{\ss{PY}}$ defined by
\[
	\mu_{\ss{PY}}(B)\;=\;\int_{B} h_\Delta\,d\mu
\]
for every Borel set $B\subseteq{\Sigma_{A}^{+}}$ is called the
\emph{Pianigiani-Yorke measure} of the subsystem
$(\Sigma_{\Delta},{\S})$.  This measure is fully supported
on~$\cZ_{\Delta}$.  

The following is a result from~\cite{CCC}.

\begin{proposition}
Let $h_\Delta$ be the function in Proposition~\ref{prop:cms}.  We have
\[
	\lim_{n\to\infty} e^{-n P_\Delta}\,\mu(\Delta_{n})\;=\; \int
	h_\Delta\,d\mu\; =\; \mu_{\ss{PY}}({{\Sigma_{A}^{+}}})\;.
\]
\end{proposition}

\section{The Periodic Case}
Let us consider the case when $\Sigma_{\Delta}$ is irreducible but
periodic with period $m>1$.  In this case there exists a decomposition
of $\Delta=\Delta_{0}\cup\cdots\cup\Delta_{m-1}$ with the property
that if $i\in\Delta_{s}$, $j\in\Delta_{s'}$ are given such that
$A(i,j)=1$ then necessarily $s'=s+1\zmod{m}$.  This induces a disjoint
partition of $\Sigma_{\Delta}=\Omega_{0} \cup\cdots\cup \Omega_{m-1}$
such that $\S(\Omega_{s})=\Omega_{s+1\zmod{m}}$, which is the
so-called cyclically moving partition of $\Sigma_{\Delta}$.  From the
classical Ruelle-Perron-Frobenius theory (cf.~\cite{Rue}, \cite{PP}),
we know that there exist non-negative H\"older continuous functions
$w_{0},\ldots, w_{m-1}$ defined on $\Sigma_{\Delta}$, and mutually
singular probability measures%
\footnote{In fact, $\nu_{j}$ is the equilibrium state of the potential 
$S_{m}(\phi)=\sum_{i=0}^{m-1}\,\phi\scirc\S^{i}$ restricted to 
$(\Omega_{j},\S^{m})$ and $\nu_{j+1\zmod{m}}=\nu_{j}\scirc\S^{-1}$.} %
$\nu_{0},\ldots,\nu_{m-1}$ with 
$\nu_{j}$ supported on $\Omega_{j}$ satisfying
\[
	\L_{\phi_{\Delta}}(w_{j}) \;=\; e^{P_{\Delta} } \,w_{j+1\zmod{m}} 
	\;,
\]
and $\mbox{supp}(w_{j})=\Omega_{j}$ for $j=0,\ldots,m-1$.  Moreover, 
$w_{\Delta}=\sum_{j=0}^{m-1}\,w_{j}$ is strictly positive on 
$\Sigma_{\Delta}$ and
\[
	\Big\Vert e^{-nP_{\Delta}}\L_{\phi_{\Delta}}^n(\psi) - 
	\sum_{j=0}^{m-1}\, w_{j+n\zmod{m}} \int_{\Omega_{j}} \psi \: 
	d\nu_{j} \Big\Vert_{{\Sigma_{\Delta}}} \; \underset{n \to 
	\infty}{\longrightarrow}\; 0 \;,
\]
for all $\psi \in C({\Sigma_{\Delta}})$.  Putting $\psi=\1$ and 
replacing $n$ by $nm$ in the above expression, we note that 
$w_{\Delta}$ could have been defined uniquely by
\[
	w_{\Delta} \;=\; \lim_{n\to\infty}\, e^{-nmP_{\Delta}}\, 
	\L_{\phi_{\Delta}}^{nm}(\1) \;,
\]
and then $\L_{\phi_{\Delta}}(w_{\Delta})=e^{P_{\Delta}}\,w_{\Delta}$.  
Now we transfer these results from $\Sigma_{\Delta}$ to the whole 
space $\Sigma_{A}^{+}$.  Let $\cZ_{\Delta}$ be defined 
by~(\ref{eq:czdef}).  Since again $\cZ_{\Delta}$ is a non-empty finite 
union of cylinders of~$\Sigma_{A}^{+}$, we have $\mu(\cZ_{\Delta})>0$.  

Define the constants $d_{j}$~by
\[
	d_{j} \;=\; \int_{\Omega_{j+1\zmod{m}}} \L_{\Delta}(\1)\, 
	d\nu_{j+1\zmod{m}} \;,
\]
for $j=0,\ldots,m-1$, we see that $d_{j}>0$ for all $j$.
Define also the constants $\alpha_{j}(k)$ 
by $\alpha_{j}(0)=1$ and for $1\leq k\leq m-1$, 
\[
	\alpha_{j}(k) \;=\; e^{-kP_{\Delta}}\, \prod_{s=0}^{k-1} 
	d_{j+s\zmod{m}} \;,
\]
for $j=0,\ldots,m-1$.  The next is our main result.
\begin{theorem}
There exist a \emph{unique} choice of non-negative H\"older continuous 
functions $h_{0},\ldots, h_{m-1}$ defined on the \emph{whole space} 
$\Sigma_{A}^{+}$ satisfying
\[
	\L_{\Delta}(h_{j}) \;=\; e^{P_{\Delta} } \,h_{j+1\zmod{m}} \;,
\]
and $h_{j}|_{\Omega_{j}}\equiv w_{j}$ for $j=0,\ldots,m-1$.  The 
function $h_{\Delta}=\sum_{j=0}^{m-1}\,h_{j}$ is strictly positive on 
$\cZ_{\Delta}$ and it is zero on the complement $\cZ_{\Delta}^{c}$.  
Moreover,
\begin{equation}
	\Big\Vert e^{-nP_{\Delta}}\L_{\Delta}^n(\psi) - 
	\sum_{j=0}^{m-1}\, \alpha_{j}(n\zmod{m})\,
	h_{j+n\zmod{m}} \int_{\Omega_{j}} \psi \: 
	d\nu_{j} \Big\Vert_{{\Sigma_{A}^{+}}} \; \underset{n \to 
	\infty}{\longrightarrow}\; 0 \;,
	\label{eq:final}
\end{equation}
for all $\psi \in \mathcal{C}({\Sigma_{A}^{+}})$.
\label{thm:cmsper}
\end{theorem}
In particular taking $\psi=1$ and integrating the above expression
with respect to $\mu$ we obtain
\begin{corollary}
Let $\Sigma_{\Delta}$ be an irreducible and periodic subsystem of
finite type with period $m$.  The sets $\Delta_{n}$ have the following
asymptotic behaviour:
\[
	\left| e^{-nP_{\Delta}}\mu(\Delta_{n}) -
	\sum_{j=0}^{m-1}\, \alpha_{j}(n\zmod{m})\,
		\int h_{j+n\zmod{m}}\,d\mu \: 
		\right| \; \underset{n \to 
		\infty}{\longrightarrow}\; 0\;.
\]
In particular, for each $k=0,1,\ldots,m-1$ we have
\[
	\lim_{n\to\infty} e^{-(k+nm)P_{\Delta}}\mu(\Delta_{k+nm}) \;=\;
	\sum_{j=0}^{m-1}\, \alpha_{j}(k)\,
			\int h_{j+k\zmod{m}}\,d\mu \;.
\]
\end{corollary}

In Section~\ref{sec:ex} we give an explicit example where the above
numbers differ for different choices of $k$, which shows that
$e^{-nP_{\Delta}}\mu(\Delta_{n})$ does not converge in general as
$n\to\infty$ when $\Sigma_{\Delta}$ is periodic. %

\section{Proof of Theorem~\ref{thm:cmsper}}
We recall some results from~\cite{CMS}.  Let
$\mathcal{C}_{p}^{+}(\Sigma_{A}^{+})$ be the set of strictly positive
$p$-cylindrical functions (i.e.~a function depending only on the first
$p$ coordinates of the point).  Let $0<\theta<1$ be the H\"older
exponent of the potential $\phi$.  Let $\cZ_{\Delta}$ be defined as
in~(\ref{eq:czdef}).  Let $\mathcal{C}({\cZ_{\Delta}})$ denote the set
of continuous functions defined on $\cZ_{\Delta}$.  The proof of the
following Lemma can be obtained from
Appendix C in~\cite{CCC}. %
\begin{lemma}
For any $f\in\cup_{p\geq 1}\,\mathcal{C}_{p}^{+}(\Sigma_{A}^{+})$, we 
have
\begin{itemize}
\item[(i)] $\{e^{-nP_{\Delta}}\,\L_{\Delta}^{n}f\}_{n\geq 0}$ is a 
Cauchy sequence in $\mathcal{C}(\Sigma_{A}^{+})\:$;
\vspace{5pt}
\item[(ii)] $h_{\Delta} = \d\lim_{n\to\infty}\, 
\dfrac{e^{-nP_{\Delta}}\,\L_{\Delta}^{n}f}{\int f\,d\mu_{\Delta}}$ 
does not depend on the function $f\in\cup_{p\geq 
1}\,\mathcal{C}_{p}^{+}(\Sigma_{A}^{+})$ and it satisfies 
\[
	\L_{\Delta}(h_{\Delta}) \;=\; e^{P_{\Delta}}\, h_{\Delta} \;.
\]
\end{itemize}
\label{lem:cms3}
\end{lemma}

Although not explicitly mentioned in~\cite{CMS}, the function
$h_{\Delta}$ is a H\"older continuous function with the same H\"older
exponent of the potential $\phi$.  We also note that, since $\mu$ is
fixed by the dual operator of $\L$ we have
\begin{equation}
	\int_{\Delta_{n}} f\cdot g\scirc \S^{n}\, d\mu\;=\; \int 
	\L_{\phi}^{n}(\bigchi_{\Delta_{n}}\cdot f\cdot g\scirc \S^{n})\, d\mu 
	\;=\; \int g \cdot \L_{\Delta}^{n}(f) \,d\mu \;.
	\label{eq:ldelta}
\end{equation}

Now consider the case when $\Sigma_{\Delta}$ is irreducible but
periodic with period $m>1$.  Consider the decomposition of
$\Delta=\Delta_{0}\cup\cdots\cup\Delta_{m-1}$ with the property that
if $i\in\Delta_{s}$, $j\in\Delta_{s'}$ are given such that $A(i,j)=1$
then necessarily $s'=s+1\zmod{m}$.  Consider also the corresponding
cyclically moving partition $\Sigma_{\Delta}=\Omega_{0} \cup\cdots\cup
\Omega_{m-1}$ such that $\S(\Omega_{s})=\Omega_{s+1\zmod{m}}$,
i.e.~defining
\[
	\Omega_{j} \;=\; \{x\in\Sigma_{\Delta}\colon\: 
	x_{0}\in\Delta_{j}\} \;,
\]
for $j=0,\ldots,m-1$.  Let $\V^{(m)}$ be the sub-alphabet of $\V^{m}$ 
defined by
\[
	\V^{(m)} \;=\; \{(i_{0},\ldots,i_{m-1})\in \V^{m} \colon\: 
	i_{0}\to \cdots \to i_{m-1}\;\; \mbox{in}\;\; \mathcal{G} \} \;,
\]
where $\mathcal{G}$ is the defining graph of~$\Sigma_{A}^{+}$. 
Consider the transition matrix $A^{(m)}$ indexed by $\V^{(m)}\times 
\V^{(m)}$ given by
\[
	A^{(m)}\big( (i_{0},\ldots,i_{m-1}), (j_{0},\ldots,j_{m-1}) \big) 
	\;=\; 1 \quad \mbox{if} \;\; i_{m-1}\to j_{0} \;\; \mbox{in}\;\;
	\mathcal{G} \;.
\]
Using the identification
\begin{equation}
	\big( (x_{0},\ldots,x_{m-1}), (x_{m},\ldots,x_{2m-1}), \ldots 
	\big) \;\longleftrightarrow \; (x_{0},x_{1},x_{2},\ldots) \;,
	\label{eq:id}
\end{equation}
the shift transformation $\S_{m}$ on $\Sigma_{A^{(m)}}^{+}$ is 
naturally topologically conjugate to $\S^{m}$ on $\Sigma_{A}^{+}$.  In 
what follows we will abuse the notation and freely identify these 
transformations and spaces.

The normalised potential $\phi$ on 
$\Sigma_{A}^{+}$ naturally defines a potential $\phi^{(m)}$ on 
$\Sigma_{A^{(m)}}^{+}$ by 
$\phi^{(m)}=S_{m}(\phi)=\sum_{j=0}^{m-1}\,\phi\scirc\S^{j}$.  Note 
that $\phi^{(m)}$ is a normalised potential for $\S_{m}$, i.e.~  
\[
	\L_{\phi^{(m)}}(\1)(x) \;=\; \sum_{y\in \S_{m}^{-1}(x)} \, 
	e^{\phi^{(m)}(y)} \;=\; 1 \;,
\]
for all $x\in\Sigma_{A^{(m)}}^{+}$.  Now the sub-alphabet $\Delta$ of 
$\V$ defines a sub-alphabet $\Delta^{(m)}$ of $\V^{(m)}$ by
\[
	\Delta^{(m)} \;=\; \{(i_{0},\ldots,i_{m-1})\in \V^{(m)}\colon \: 
	i_{s}\in\Delta\,, \; \mbox{for}\; s=0,\ldots, m-1 \}\;,
\]
and this sub-alphabet can be further decomposed into
\begin{equation}
	\Delta^{(m)}_{j} \;=\; \{(i_{0},\ldots,i_{m-1})\in 
	\Delta^{(m)}\colon \: i_{s}\in\Delta_{s+j\zmod{m}}\,, \; 
	\mbox{for} \; s=0,\ldots,m-1 \} \;,
	\label{eq:dmj}
\end{equation}
for $j=0,\ldots,m-1$.  The important fact is that for fixed $j$, 
$\Sigma_{j}^{(m)}$ the subsystem of~$\Sigma_{A^{(m)}}^{+}$ obtained by 
taking transitions through $\Delta^{(m)}_{j}$ is irreducible and 
aperiodic in its alphabet $\Delta^{(m)}_{j}$.  Hence the main result 
of~\cite{CMS} applies and we define a H\"older continuous function 
$h_{j}$ by
\begin{equation}
	h_{j} \;=\; \lim_{n\to\infty}\, 
	e^{-nP(\Delta^{(m)}_{j})}\,\L_{\Delta^{(m)}_{j}}^{n}(\1) \;,
	\label{eq:hjdef}
\end{equation}
where $P(\Delta^{(m)}_{j})$ is the pressure of the restriction of 
$\phi^{(m)}$ to the subsystem $\Sigma^{(m)}_{j}$ (hence 
$P(\Delta^{(m)}_{j})=mP_{\Delta}$ for all $j$), and 
$\L_{\Delta^{(m)}_{j}}(\psi)= 
\L_{\phi^{(m)}}(\psi\cdot\bigchi_{\Delta^{(m)}_{j}})$ with respect to 
the shift $\S_{m}$.  From Lemma~\ref{lem:cms3}~(ii) extended to 
continuous functions we obtain
\begin{equation}
	\lim_{n\to\infty} \, e^{-nmP_{\Delta}}\, 
	\L_{\Delta^{(m)}_{j}}^{n}(\psi) \;=\; h_{j} 
	\int_{\Sigma^{(m)}_{j}} \psi\,d\nu_{j} \;,
	\label{eq:hjint}
\end{equation}
for every $\psi\in\mathcal{C}(\Sigma_{A^{(m)}}^{+})$ (with the limit 
being uniform), where $\nu_{j}$ is the unique equilibrium state of 
$\phi^{(m)}$ restricted to the subsystem $\Sigma^{(m)}_{j}$ with 
respect to the shift $\S_{m}$.

In view of the identification~(\ref{eq:id}) the function $h_{j}$ 
defines a function on~$\Sigma_{A}^{+}$ in a natural way and $\nu_{j}$ 
becomes a probability measure on~$\Sigma_{A}^{+}$ fully supported 
on~$\Omega_{j}$.  Note that $h_{j}$ is then strictly positive on
\[
	\cZ_{\Delta_{j}} \;=\; \{x\in\Sigma_{A}^{+}\colon\: 
	\exists\,b\in\Delta_{j-1\zmod{m}}\,, \;A(b,x_{0})=1 \} \;,
\]
and it is zero on the complement $\cZ_{\Delta_{j}}^{c}$.  Note also 
that for each $j$, $\cZ_{\Delta_{j}}$ is a non-empty finite union of 
cylinders of~$\Sigma_{A}^{+}$, therefore in particular, 
$\mu(\cZ_{\Delta_{j}})>0$. Applying $\L_{\Delta}$ as
\[
\begin{split}
	\L_{\Delta}(\psi)\big( (x_{0},\ldots,x_{m-1}), 
	(x_{m},&\ldots,x_{2m-1}), \ldots \big)\; \\ %
	=\; \sum_{\{i\in\Delta\colon A(i,x_{0})=1\}} &e^{\phi^{(m)} ( 
	(i,x_{0},\ldots,x_{m-2}), (x_{m-1},\ldots,x_{2m-2}), \ldots )} 
	\;\times \\ % 
	&\psi\big( (i,x_{0},\ldots,x_{m-2}), (x_{m-1},\ldots,x_{2m-2}), 
	\ldots \big) \;,
\end{split}
\]
we conclude that $\L_{\Delta}\scirc\L_{\Delta^{(m)}_{j}}^{n} = 
\L_{\Delta^{(m)}_{j+1\zmod{m}}}^{n}\!\!\!\scirc\L_{\Delta}$.  This 
implies that for all $k\geq 1$ we have 
$\L_{\Delta}^{k}\scirc\L_{\Delta^{(m)}_{j}}^{n} = 
\L_{\Delta^{(m)}_{j+k\zmod{m}}}^{n}\!\!\!\scirc\L_{\Delta}^{k}$.  
Putting $\psi=\1$ in~(\ref{eq:hjint}) we obtain, for fixed $k\geq 1$ 
and fixed $0\leq j<m$,
\begin{equation}
\begin{split}
	\L_{\Delta}^{k}(h_{j}) \;&=\; \lim_{n\to\infty} \, 
	e^{-nmP_{\Delta}}\, \L_{\Delta}^{k}(\L_{\Delta^{(m)}_{j}}^{n}(\1)) 
	\;=\; \lim_{n\to\infty} \, e^{-nmP_{\Delta}}\, 
	\L_{\Delta^{(m)}_{j+k\zmod{m}}}^{n}(\L_{\Delta}^{k}(\1)) \;\\ %
	&=\; h_{j+k\zmod{m}} \int_{\Sigma^{(m)}_{j+k\zmod{m}}} 
	\!\!\!\L_{\Delta}^{k}(\1)\,d\nu_{j+k\zmod{m}} \;.
\end{split}
\label{eq:lkhj}
\end{equation}
Therefore defining the constants $d_{j}$~by
\[
	d_{j} \;=\; \int_{\Omega_{j+1\zmod{m}}} \L_{\Delta}(\1)\, 
	d\nu_{j+1\zmod{m}} \;,
\]
for $j=0,\ldots,m-1$, we see that $d_{j}>0$ for all $j$ 
and~by~(\ref{eq:lkhj}) we have
\[
	\L_{\Delta}(h_{j}) \;=\; d_{j}\,h_{j+1\zmod{m}} \;,
\]
for all $j$.  Since 
$\L_{\Delta}^{m}(h_{j})=e^{mP_{\Delta}}\,h_{j}$ for each $j$, 
by~(\ref{eq:lkhj}) we also see that 
\[
	\prod_{j=0}^{m-1}\,d_{j} \;=\; e^{mP_{\Delta}} \;.
\]
At the end of this appendix we give an example where in general one 
has $d_{j}$ not necessarily equal to $e^{P_{\Delta}}$.
%  (this is already 
% present in the periodic case of classical Ruelle-Perron-Frobenius 
% theory, see \cite{Rue} or \cite{PP}).

The function $h_{\Delta}=\sum_{j=0}^{m-1} h_{j}$ is strictly positive 
on~$\cZ_{\Delta}$ and it is zero on the complement $\cZ_{\Delta}^{c}$.  
(Note that $\cZ_{\Delta}=\cup_{j=0}^{m-1}\cZ_{\Delta_{j}}$ and this 
union is not in general a disjoint union, see example below.)  
The function $h_{\Delta}$ also satisfies
\[
	\L_{\Delta}^{m}(h_{\Delta}) \;=\; e^{mP_{\Delta}}\,h_{\Delta} \;.
\]
Now, from the fact that
\[
	\L_{\Delta^{(m)}}(\psi) \;=\; \sum_{j=0}^{m-1}\, 
	\L_{\Delta^{(m)}_{j}}(\psi) \;,
\]
and $\L_{\Delta^{(m)}_{j}}\scirc\L_{\Delta^{(m)}_{j'}}=0$ if 
$j\not=j'$, we see that
\[
	\L_{\Delta^{(m)}}^{n}(\psi) \;=\; \sum_{j=0}^{m-1}\, 
	\L_{\Delta^{(m)}_{j}}^{n}(\psi) \;.
\]
Using~(\ref{eq:lkhj}) we have, for fixed $1\leq k<m$,
\begin{equation}
\begin{split}
	e^{-(nm+k)P_{\Delta}}\,\L_{\Delta}^{nm+k}(\psi) \;&=\; 
	e^{-(nm+k)P_{\Delta}}\,\L_{\Delta}^{k}\big( 
	\L_{\Delta^{(m)}}^{n}(\psi) \big) \\ %
	&=\; \sum_{j=0}^{m-1}\, e^{-kP_{\Delta}}\, \L_{\Delta}^{k}\big( 
	e^{-nmP_{\Delta}} \L_{\Delta^{(m)}_{j}}^{n}(\psi) \big) \\ %
	&=\; \sum_{j=0}^{m-1}\, e^{-kP_{\Delta}}\, \L_{\Delta}^{k}(h_{j}) 
	\,\int_{\Omega_{j}} \psi\,d\nu_{j} + o(1) \\ %
	&=\; \sum_{j=0}^{m-1}\, \Big( e^{-kP_{\Delta}}\, \prod_{s=0}^{k-1} 
	d_{j+s\zmod{m}} \Big) \,h_{j+k\zmod{m}} \,\int_{\Omega_{j}} 
	\psi\,d\nu_{j} + o(1) \;,
\end{split}	
\label{eq:bdper}
\end{equation}
where we have used Lemma~\ref{lem:cms3}~(ii) for continuous functions 
and $o(1)$ is with respect to $n$.  Define the constants $\alpha_{j}(k)$ 
by $\alpha_{j}(0)=1$ and for $1\leq k\leq m-1$, 
\[
	\alpha_{j}(k) \;=\; e^{-kP_{\Delta}}\, \prod_{s=0}^{k-1} 
	d_{j+s\zmod{m}} \;,
\]
for $j=0,\ldots,m-1$.  Therefore from~(\ref{eq:bdper}) we finally 
obtain
\[
	\Big\Vert e^{-nP_{\Delta}}\L_{\Delta}^n(\psi) - 
	\sum_{j=0}^{m-1}\, \alpha_{j}(n\zmod{m})\,
	h_{j+n\zmod{m}} \int_{\Omega_{j}} \psi \: 
	d\nu_{j} \Big\Vert_{{\Sigma_{A}^{+}}} \; \underset{n \to 
	\infty}{\longrightarrow}\; 0 \;,
\]
for all $\psi \in \mathcal{C}({\Sigma_{A}^{+}})$, which concludes the 
proof of Theorem~\ref{thm:cmsper}.

\section{Illustrative Example}
\label{sec:ex}
In this section we give an example to illustrate the computations made
in the previous section.
\begin{example}
Let $\V=\{1,2,3\}$ and consider the matrix $A$ given by
\[
	A \;=\; 
	\left( \begin{array}{ccc}
	0 & 1 & 1 \\
	1 & 0 & 1 \\
	1 & 1 & 1
	\end{array}
	\right)\;.
\]
Let $\phi$ be any normalised H\"older continuous potential on 
$\Sigma_{A}^{+}$, i.e.~assume that
\[
	\L_{\phi}(\1)(x) \;=\; \sum_{\{i\in\V\colon A(i,x_{0})=1\}} 
	e^{\phi(ix)} \;=\; 1 \;,
\]
for all $x\in\Sigma_{A}^{+}$.  Take $\Delta=\{1,2\}$.  Then 
$\Sigma_{\Delta}$ is the periodic orbit $\{(1,2,1,2,\ldots), 
(2,1,\\ 2,1,\ldots)\}$.  Put $\phi_{\Delta}(1,2,1,2,\ldots)=p$ and 
$\phi_{\Delta}(2,1,2,1,\ldots)=q$, and assume $p\not=q$.  There is 
only one invariant measure for the restriction of the shift $\S$ on 
$\Sigma_{\Delta}$, namely
\[
	\mu_{\Delta} \;=\; \dfrac{1}{2}(\delta_{1}+\delta_{2})\;,
\]
where $\delta_{1}$ is Dirac measure at the point $(1,2,1,2,\ldots)$ 
and $\delta_{2}$ is Dirac measure at the point $(2,1,2,1,\ldots)$.  
Since the shift entropy of $\mu_{\Delta}$ is zero, the restricted 
pressure $P_{\Delta}$ is then given by
\[
	P_{\Delta} \;=\; \int \phi_{\Delta}\,d\mu_{\Delta} \;=\; 
	\dfrac{1}{2}(p+q)\;.
\]
Notice now that $\phi^{(2)}=\phi + \phi\scirc\S$ when restricted 
to~$\Sigma_{\Delta}$ is constant with value $p+q$.  The pressure of 
the restriction of $\phi^{(2)}$ to $\Sigma_{\Delta}$ with respect to 
$\S^{2}$ is then given by~$p +q =2P_{\Delta}$.  The set $\Delta$ is 
further decomposed into $\Delta_{i-1}=\{i\}$, for $i=1,2$, giving the 
cyclically moving partition 
$\Sigma_{\Delta}=\Omega_{0}\cup\Omega_{1}$, where 
$\Omega_{0}=\{(1,2,1,2,\ldots)\}$ and 
$\Omega_{1}=\{(2,1,2,1,\ldots)\}$.  Note then that $T^{2}$ restricted 
to~$\Sigma_{\Delta}$ consists of two fixed points.  This implies that 
$\nu_{0}=\delta_{1}$ and $\nu_{1}=\delta_{2}$, where $\nu_{i}$ is the 
equilibrium state of $\phi^{(2)}$ restricted to $\Omega_{i}$ with 
respect to~$\S^{2}$.  Applying~\cite{CMS} in the case of an aperiodic 
subsystem consisting of a fixed point for $\S^{2}$ we have 
from~(\ref{eq:hjdef}), where $\Delta^{(m)}_{j}$ is defined 
by~(\ref{eq:dmj}) and $m=2$,
\[
	h_{j} \;=\; \lim_{n\to\infty}\, e^{-2nP_{\Delta}}\, 
	\L_{\Delta^{(2)}_{j}}^{n}(\1) \;,
\]
for $j=0,1$.  Interpreting this we conclude that for 
$x\in\cZ_{\Delta_{0}}=C[1]_{0}\cup C[3]_{0}$ we have
\[
\begin{split}
	h_{0}(x) \;&=\; \lim_{n\to\infty}\, \exp\big\{ S_{2n}(\phi)\big( 
	\underbrace{1,2,1,2,\ldots 1,2}_{2n},x_{0},x_{1},\ldots\big) - 
	n(p+q) \big\} \\ %
	&=\; \lim_{n\to\infty}\, \exp\big\{ S_{2n}(\phi)\big( 
	\underbrace{1,2,1,2,\ldots 1,2}_{2n},x_{0},x_{1},\ldots\big) - 
	S_{2n}(\phi)\big(1,2,1,2,\ldots \big)\big\}
	\;,
\end{split}
\]
where $S_{k}(\phi)$ denotes 
$\phi+\phi\scirc\S+\cdots+\phi\scirc\S^{k-1}$, and $h_{0}$ is zero on 
the complement $\cZ_{\Delta_{0}}^{c}=C[2]_{0}$.  Also 
if~$x\in\cZ_{\Delta_{1}}=C[2]_{0}\cup C[3]_{0}$ then
\[
\begin{split}
	h_{1}(x) \;&=\; \lim_{n\to\infty}\, \exp\big\{ S_{2n}(\phi)\big( 
	\underbrace{2,1,2,1,\ldots 2,1}_{2n} ,x_{0},x_{1},\ldots\big) - 
	n(p+q) \big\} \\ %
	&=\; \lim_{n\to\infty}\, \exp\big\{ S_{2n}(\phi)\big( 
	\underbrace{2,1,2,1,\ldots 2,1}_{2n} ,x_{0},x_{1},\ldots\big) - 
	S_{2n}(\phi)\big( 2,1,2,1,\ldots\big) \big\}
	\;,
\end{split}
\]
and $h_{1}$ is zero on the 
complement~$\cZ_{\Delta_{1}}^{c}=C[1]_{0}$. (Note that $h_{0}$ and 
$h_{1}$ are both strictly positive on the cylinder~$C[3]_{0}$.) Now 
we compute the constants $d_{j}$, for $j=0,1$. We have
\[
\begin{split}
	d_{0} \;&=\; \int_{\Omega_{1}} \L_{\Delta}(\1)\, 
	d\nu_{1} \;=\; e^{\phi(1,2,1,2,\ldots)} \;=\; e^{p} \;, 
	\qquad\mbox{and} \\ %
	d_{1} \;&=\; \int_{\Omega_{0}} \L_{\Delta}(\1)\, 
	d\nu_{0} \;=\; e^{\phi(2,1,2,1,\ldots)} \;=\; e^{q} \;.
\end{split}
\]
This provides an example where 
$d_{j}\not=e^{P_{\Delta}}=e^{\frac{1}{2}(p+q)}$, since we are assuming 
$p\not= q$. One can see directly that $h_{0}$ and $h_{1}$ satisfy
\[
	\L_{\Delta}(h_{0}) \;=\; d_{0}\,h_{1} \;=\; e^{p}\,h_{1} 
	\quad\mbox{and}\quad 
	\L_{\Delta}(h_{1}) \;=\; d_{1}\,h_{0} \;=\; e^{q}\,h_{0} \;.
\]
The function $h_{\Delta}=h_{0}+h_{1}$ satisfies 
$\L_{\Delta}^{2}(h_{\Delta})=e^{2P_{\Delta}}\,h_{\Delta}$, and in the 
case of this example it is fully supported on~$\Sigma_{A}^{+}$.  Now 
we compute the constants $\alpha_{j}(k)$ for $j,k=0,1$.  We have 
$\alpha_{j}(0)=1$ for $j=0,1$,
\[
\begin{split}
	\alpha_{0}(1) \;&=\; e^{-P_{\Delta}}\, d_{0} \;=\; 
	e^{-\frac{1}{2}(p+q)}\, e^{p} \;=\; e^{\frac{1}{2}(p-q)} \;,
	\qquad\mbox{and} \\ %
	\alpha_{1}(1) \;&=\; e^{-P_{\Delta}}\, d_{1} \;=\; 
	e^{-\frac{1}{2}(p+q)}\, e^{q} \;=\; e^{\frac{1}{2}(q-p)} 
	\;.
\end{split}
\]
From~(\ref{eq:final}) we conclude that
\[
\begin{split}
	\Big\Vert e^{-nP_{\Delta}}&\L_{\Delta}^n(\psi) - 
	\Big( \alpha_{0}(n\zmod{m})\, h_{n\zmod{m}} 
	\,\psi(1,2,1,2,\ldots) \\ % 
	&+ \alpha_{1}(n\zmod{m})\, h_{n+1\zmod{m}} 
	\,\psi(2,1,2,1,\ldots) \Big) \Big\Vert_{{\Sigma_{A}^{+}}} \; 
	\underset{n \to \infty}{\longrightarrow}\; 0 \;,
\end{split}
\]
for all $\psi \in \mathcal{C}({\Sigma_{A}^{+}})$. An interesting fact 
is that, putting $f=g=\1$ in~(\ref{eq:ldelta}) and putting $\psi=\1$ 
in the above expression we have
\[
\begin{split}
	\mu(\Delta_{n}) \;&=\; \int_{\Delta_{n}}\,d\mu \;=\; 
	\int\L_{\Delta}^{n}(\1)\,d\mu \\ %
	&=\; e^{nP_{\Delta}}\,\Big( \alpha_{0}(n\zmod{m})\, \int 
	h_{n\zmod{m}}\,d\mu \\ %
	&+ \alpha_{1}(n\zmod{m})\, \int h_{n+1\zmod{m}}\,d\mu
	\Big) + o(e^{nP_{\Delta}})\;.
\end{split}
\]
Therefore
\[
\begin{split}
	\lim_{n\to\infty}\, e^{-2nP_{\Delta}}\,\mu(\Delta_{2n}) \;&=\; 
	\alpha_{0}(0)\, \int 
	h_{0}\,d\mu + \alpha_{1}(0)\, \int h_{1}\,d\mu \\ %
	&=\; \int (h_{0}+h_{1})\,d\mu \;=\; \int h_{\Delta}\,d\mu \;,
\end{split}
\]
but
\[
\begin{split}
	\lim_{n\to\infty}\, e^{-(2n+1)P_{\Delta}}\,\mu(\Delta_{2n+1}) 
	\;&=\; \alpha_{0}(1)\, \int h_{1}\,d\mu + \alpha_{1}(1)\, \int 
	h_{0}\,d\mu \\ % 
	&=\; e^{\frac{1}{2}(q-p)}\int h_{0}\,d\mu +
	e^{\frac{1}{2}(p-q)}\int h_{1}\,d\mu \;.
\end{split}
\]
The latter is not in general equal to $\int h_{\Delta}\,d\mu$ if 
$p\not=q$ (see explicit example below).  Therefore 
$\lim_{n\to\infty}\, e^{-nP_{\Delta}}\,\mu(\Delta_{n})$ may not exist 
in general.  However, if $p=q$ then $\alpha_{j}(k)=1$ for all $j,k$ 
and then the limit is given by
\[
	\lim_{n\to\infty}\, e^{-nP_{\Delta}}\,\mu(\Delta_{n}) \;=\; \int 
	h_{\Delta}\,d\mu \;.
\]
Note also that even when~$p\not=q$ there are choices of normalised 
potential $\phi$ such that
\[
	\lim_{n\to\infty}\, e^{-nP_{\Delta}}\,\mu(\Delta_{n}) \;=\;
	e^{\frac{1}{2}(q-p)}\int h_{0}\,d\mu + e^{\frac{1}{2}(p-q)}\int 
	h_{1}\,d\mu \;=\; \int h_{\Delta}\,d\mu\;.
\]
For explicit examples of the above remarks, take for instance $\phi$ 
defined by $\phi|_{C[1]_{0}}\equiv p$ and $\phi|_{C[2]_{0}}\equiv q$.  
Then necessarily $h_{0}$ is equal to $1$ on the cylinders $C[1]_{0}$ 
and $C[3]_{0}$, and it is equal to $0$ on $C[2]_{0}$.  Similarly, 
$h_{1}$ is equal to $1$ on $C[2]_{0}$ and $C[3]_{0}$, and it is equal 
to $0$ on $C[1]_{0}$.  Therefore
\[
\begin{split}
	\int h_{\Delta}\,d\mu \;&=\; \big( \mu(C[1]_{0}) + \mu(C[3]_{0}) 
	\big) + \big( \mu(C[2]_{0}) + \mu(C[3]_{0}) \big) \\ %
	&=\; \big(1-\mu(C[2]_{0}) \big) + \big(1-\mu(C[1]_{0}) \big) \;.
\end{split}
\]
Now the condition of $\phi$ being normalised implies that the values 
of $\phi$ on the cylinder $C[3]_{0}$ is uniquely determined.  In fact 
on this cylinder $\phi$ is the 2-step cylindrical function given by
\[
	\phi|_{C[31]_{0}} \;\equiv\; \log(1-e^{q})\,, \;\; 
	\phi|_{C[32]_{0}} \;\equiv\; \log(1-e^{p})\,, \;\; \mbox{and}\;\;
	\phi|_{C[33]_{0}} \;\equiv\; \log(1-e^{p}-e^{q})\;.
\]
Hence $\mu$ is the Markov measure defined by the stochastic matrix
\[
	P \;=\; \left(
	\begin{array}{ccc}
	0 & e^{q} & 1-e^{q} \\
	e^{p} & 0 & 1-e^{p} \\
	e^{p} & e^{q} & 1-e^{p}-e^{q}
	\end{array}
	\right) \;.
\]
This matrix has the stationary strictly positive left eigenvector 
$(p_{1},p_{2},p_{3})$ given by
\[
	(p_{1},p_{2},p_{3}) \;=\; \left(
	\dfrac{e^{p}}{1+e^{p}}\, ,\: \dfrac{e^{q}}{1+e^{q}}\,, \:
	1- \dfrac{e^{p}}{1+e^{p}} - \dfrac{e^{q}}{1+e^{q}}
	\right) \;.
\]
Therefore, $\mu(C[i]_{0})=p_{i}$ for $i=1,2,3$, which implies that
\[
	\int h_{\Delta}\,d\mu \;=\; \left(1 - \dfrac{e^{p}}{1+e^{p}} 
	\right) + \left(1 - \dfrac{e^{q}}{1+e^{q}} 
	\right) \;=\; \dfrac{2+e^{p}+e^{q}}{(1+e^{p})(1+e^{q})} \;.
\] 
Now we compare the above expression with
\[
\begin{split}
	e^{\frac{1}{2}(q-p)}\int h_{0}\,d\mu &+
	e^{\frac{1}{2}(p-q)}\int h_{1}\,d\mu \;=\; e^{\frac{1}{2}(q-p)}\, 
	(1-p_{2}) + e^{\frac{1}{2}(p-q)}\, (1-p_{1}) \\ %
	&=\; \dfrac{e^{\frac{1}{2}(q-p)}}{1+e^{q}} + 
	\dfrac{e^{\frac{1}{2}(p-q)}}{1+e^{p}} \;=\; \dfrac{ 
	e^{\frac{1}{2}(q-p)}\,(1+e^{p}) + e^{\frac{1}{2}(p-q)}\,(1+e^{q}) 
	}{(1+e^{p})(1+e^{q})} \;.
\end{split}	
\]
The two expressions coincide if and only if
\[
	2+e^{p}+e^{q} \;=\; e^{\frac{1}{2}(q-p)}\,(1+e^{p}) + 
	e^{\frac{1}{2}(p-q)}\,(1+e^{q}) \;.
\]
Introducing $a=e^{\frac{1}{2}(q-p)}$ we see that
\[
	2 +e^{p}+a^{2}\,e^{p} \;=\; a\,(1+e^{p}) + 
	a^{-1}\,(1+a^{2}\,e^{p}) \;,
\]
which is equivalent to
\[
	a\,(a^{2}-2a+1)\,e^{p} \;=\; a^{2}-2a+1 \;.
\]
Since $a\,e^{p}=e^{\frac{1}{2}(p+q)}=e^{P_{\Delta}}<1$, the above 
equality holds if and only if $a=1$ (i.e.~if $p=q$). Therefore, 
whenever $p\not=q$, we have 
\[
	\int 
	h_{\Delta}\,d\mu \;\not=\; 
	e^{\frac{1}{2}(q-p)}\int h_{0}\,d\mu + e^{\frac{1}{2}(p-q)}\int 
	h_{1}\,d\mu \;.
\]

Hence, $e^{-nP_{\Delta}}\mu(\Delta_{n})$ does not converge as
$n\to\infty$ when $p\neq q$.
\end{example}

\bibliographystyle{article}

\end{document}